# ON A PROBLEM OF RAMSEY THEORY

**Frasser C.E.**
*Ph.D. in Engineering Sciences, Odessa National Polytechnic University, Ukraine*

**Abstract**
In 1955, Greenwood and Gleason showed that the Ramsey number $R(3, 3, 3) = 17$ by constructing an edge-chromatic graph on 16 vertices in three colors with no triangles. Their technique employed finite fields. This same result was obtained later by using another technique. In this article, we examine the complete graph on 17 vertices, $K_{17}$, which can be represented as a regular polygon of 17 sides with all its diagonals. We color each edge of $K_{17}$ with one of the three colors, blue, red or yellow. The graph thus obtained is called *complete trichromatic graph* $K_{17}^{(3)}$ (the superscript determines the number of colors). A triangle contained in graph $K_{17}^{(3)}$ with edges colored with one and only one color is called *monochromatic*. It has been shown that for any coloring of the $K_{17}^{(3)}$ edges, $K_{17}^{(3)}$ contains at least one monochromatic triangle. This article examines the problem of determining the minimum number of monochromatic triangles with the same color contained in $K_{17}^{(3)}$.

In order to solve the aforementioned problem, construct first the graph $K_{16}^{(3)}$ that does not contain monochromatic triangles. An example of the construction of such a graph was described for the first time in the classic research article authored by Greenwood and Gleason [2].

Let's designate the vertices of $K_{16}$ with the letters $O$, $A_i$, $B_i$, $C_i$, $i = 1,\ldots,5$ (Fig.1). Let's color the edges $OA_1,\ldots, OA_5$ ; $OB_1,\ldots, OB_5$ ; $OC_1,\ldots, OC_5$ with the colors blue, red and yellow, respectively. Let's color subgraphs $K_5$ with vertices $A_i$, $B_i$, $C_i$, $i = 1,\ldots,5$ with the colors red and yellow, yellow and blue, blue and red. Now, let's color the edges between $A_i$ and $B_j$, $B_i$ and $C_j$, $C_i$ and $A_j$ ; $i, j = 1,\ldots,5$, using the convenient color so that no monochromatic triangles are formed. Such a coloring of edges $A_iB_j$, $i, j = 1,\ldots,5$ is obtained by joining the left and right parts of Fig. 1 so that the cylinder (Fig. 2 (a), (b)) is formed.

Using the cyclic permutation of the colors

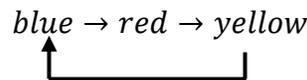

$$blue \to red \to yellow$$

we obtain the corresponding colorings for the edges $B_iC_j$, $C_iA_j$, that is, if $A_1B_2$ is red, then $B_1C_2$ is yellow and $C_1A_2$ is blue. Such a coloring guarantees the absence of monochromatic triangles with vertex $O$ and with vertices $A_i$ and $B_j$, $B_i$ and $C_j$, $C_i$ and $A_j$ ; $i, j = 1,\ldots,5$.

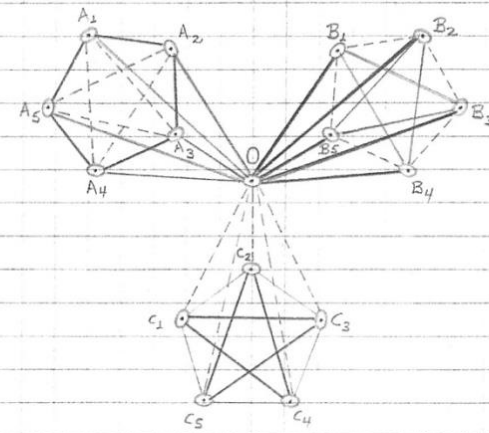

Fig. 1

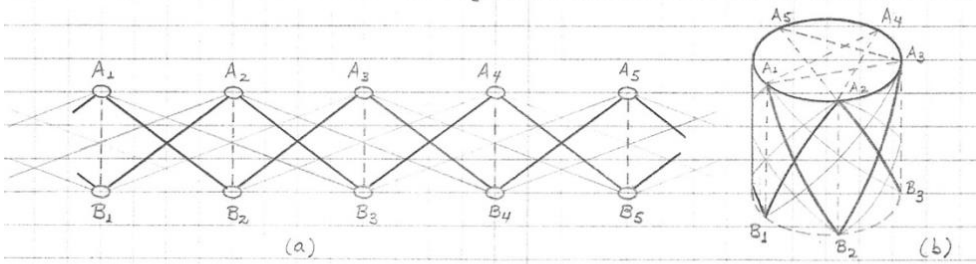

(a)  (b)

Fig. 2.

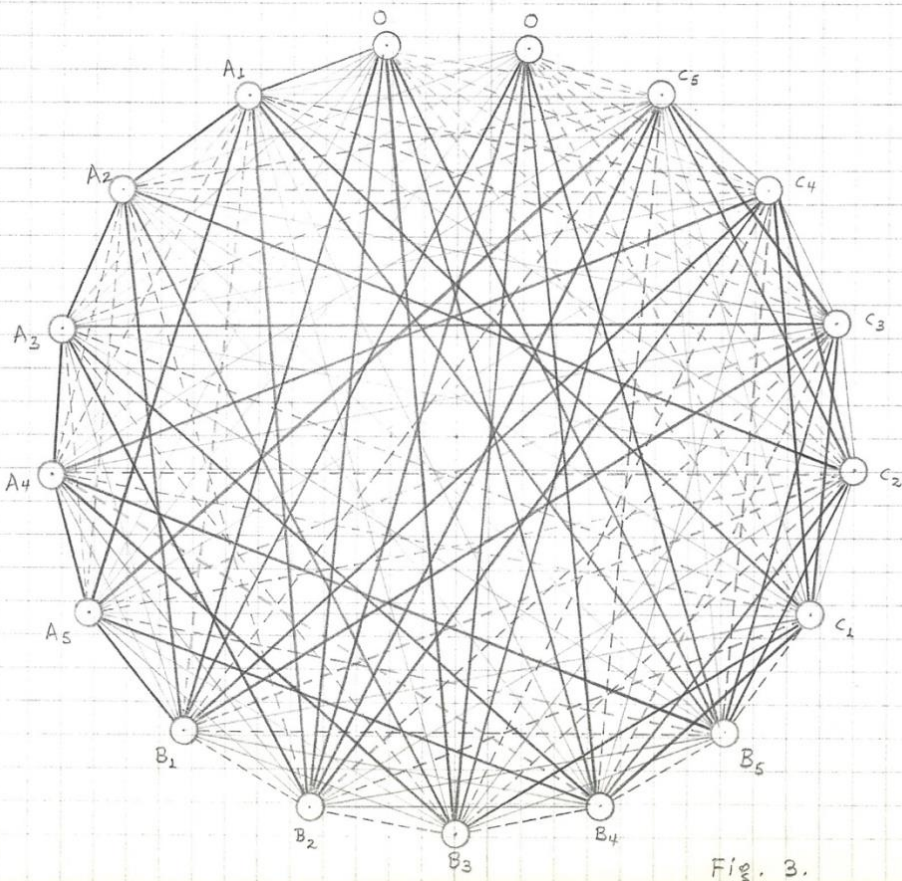

Fig. 3.

It only remains to check that there are no monochromatic triangles of the form $A_i B_j C_k$; $i, j, k = 1,\ldots,5$. It is easy to convince yourself of it by checking directly. Therefore, the constructed graph $K_{16}^{(3)}$ does not contain monochromatic triangles.

We now need to show that there exists graph $K_{17}^{(3)}$, which contains exactly 5 monochromatic triangles of only one color. To this purpose, it should result from the graph of Fig. 3 two subgraphs $K_{16}^{(3)}$ with no monochromatic triangles, in which the colors of the blue, red, and yellow edges coincide. That is only possible if when deleting vertex $O$ (including the corresponding incident edges) from each of the two subgraphs $K_{16}^{(3)}$ of the graph in Fig. 3, we obtain identical subgraphs $K_{15}^{(3)}$. As a result, we get graph $K_{17}^{(3)}$ without one of its edges, in which there is no monochromatic triangles. Adding the missing edge and coloring it blue, red or yellow, we obtain graph $K_{17}^{(3)}$ with five triangles colored blue, red or yellow, respectively. This elegant construction is due to mathematician A.D. Glukhov (Kiev, Ukraine), whose description is based on the results obtained by C.B. Beliy [1]. The original proof, omitted in this article, belongs to C.B. Beliy.